\documentclass[12pt]{article}
\usepackage{ amsmath, amsfonts, amssymb,amsthm}
\usepackage{color, graphicx}

\newtheorem{definition}{Definition}[section]
\newtheorem{lemma}[definition]{Lemma} 
\newtheorem{proposition}[definition]{Proposition}

\newtheorem{theorem}[definition]{Theorem}

\def\reg{\operatorname{reg}}
\begin{document}

\makeatletter      
\renewcommand{\ps@plain}{%
     \renewcommand{\@oddhead}{\textrm{LOWER BOUND FOR THE REGULARITY INDEX OF FAT POINTS}\hfil\textrm{\thepage}}%
     \renewcommand{\@evenhead}{\@oddhead}%
     \renewcommand{\@oddfoot}{}
     \renewcommand{\@evenfoot}{\@oddfoot}}
\makeatother     

\title{LOWER BOUND FOR  THE REGULARITY INDEX OF FAT POINTS}         
\author{Phan Van Thien}        
\date{}          
\maketitle

\pagestyle{plain}

\begin{abstract}\noindent The problem to find an upper bound for the regularity index of fat points
 has been dealt with by many authors. In this paper we give a lower bound for the regularity index of fat points.
It is useful for determining the regularity index.
\end{abstract}



\par \section{Introduction}

Let $P_1,\ldots , P_s$ be distinct points in the projective space $\Bbb P^n :=
\Bbb P^n(K)$, $K$ an algebraically closed field. Denote by $\wp_1,\ldots ,
\wp_s$ the prime ideals in the polynomial ring $R:=K[X_0,\ldots , X_n]$
corresponding to the points $P_1,\ldots , P_s$. Let $m_1,\ldots , m_s$ be
positive integers. We will denote by $Z$ the zero-scheme defined by  the ideal
$I:=\wp^{m_1}_1 \cap \cdots \cap \wp^{m_s}_s$ and call $Z$ a set of {\it fat
points} in $\Bbb P^n$.

\par The homogeneous coordinate ring of $Z$ is $R/I$. This ring is a one-dimensional Cohen-Macaulay graded ring,
$R/I=\underset{t \ge 0}{\oplus} (R/I)_t$, whose multiplicity is
$$e(R/I) = \sum^s_{i=1}\binom{m_i+n-1}{n}.$$
For every $t$, the $(R/I)_t$ is a finite dimensional vector space over $K$. The
function $H_{R/I}(t):= \text{dim}_K (R/I)_t$ strictly increases until it reaches
the multiplicity $e(R/I)$, at which it stabilizes. The {\it regularity index} of
$Z$, denote by $\reg(Z)$, is defined to be the least integer $t$ such that
$H_{R/I}(t)=e(R/I)$. It is well known that $\reg(Z)=\reg(R/I)$, the
Castelnuovo-Mumford regularity of $R/I$. Hence we will also denote $\reg(Z)$ by
$\reg(R/I)$.
\smallskip

The problem to exactly determine the regularity index $\reg(Z)$ is fairly
difficult. So, instead of determining $\reg(Z)$, one tries to find an upper
bound for it. The problem to find an upper bound for $\reg(Z)$ has been dealt
with by many authors (see \cite{Ba1}-\cite{Th6}).

In this paper we give a lower bound for the regularity index of fat points. The
lower bound and upper bound are useful tools for determining the regularity
index.

The algebraic method used in this paper as well as in \cite{CTV}.

\section{Preliminaries}

From now on, we say a $j$-plane, i.e. a linear $j$-space. We identify a
hyperplane as the linear form defining it.
\par

We will use the following lemmas which have been proved in \cite{CTV}.

\begin{lemma}\label{l21} \cite[Lemma 1]{CTV} Let $P_1,\ldots , P_r, P$ be
distinct points in $\Bbb P^n$ and let $\wp$ be the defining ideal of $P$. If
$m_1,\ldots , m_r$ and $a$ are positive integers, $J := \wp^{m_1}_1\cap
\cdots\cap  \wp^{m_r}_r$, and $I = J \cap \wp^a$, then $$\reg(R/I) =
\max\left\{a-1, \reg(R/J), \reg(R/(J+\wp^a)) \right\}.$$ \end{lemma}

\par
\begin{lemma}\label{l22} \cite[Lemma 3]{CTV}
 Let $P_1,\ldots , P_r$ be distinct
points in $\Bbb P^n$ and $a, m_1,\ldots , m_r$ positive integers. Put $J =
\wp^{m_1}_1\cap \cdots\cap  \wp^{m_r}_r$ and $\wp = (X_1,\ldots , X_n)$. Then
$$\reg(R/(J+\wp^a)) \le b$$ if and only if $X^{b-i}_0M \in J+\wp^{i+1}$ for
every monomial $M$ of degree $i$ in $X_1,\ldots , X_n$, $i = 0,\ldots , a-1$.
\end{lemma}  \par

Suppose that we can find $t$ hyperplanes $H_1, \ldots, H_t$ avoiding $P$ such
that $H_1 \cdots H_t  M \in J$ for every monomial $M$ of degree $i$ in
$X_1,\ldots , X_n$, $i = 0,\ldots , a-1$. Since we can write $H_j=X_0+G_j$ for
some linear form $G_j\in \wp$ for $j=1, \ldots, t$, we get $X_0^{t}M \in
J+\wp^{i+1}$. Therefore, we have the following lemma:

\begin{lemma}\label{l23} Assume that $H_1, \ldots, H_t$ are hyperplanes
avoiding $P$ such that $H_1 \cdots H_t  M \in J$ for every monomial $M$ of
degree $i$ in $X_1,\ldots , X_n$, $i = 0,\ldots , a-1$. If
$$\delta \ge \max\{t+i| i=1, \ldots, a-1\}$$
then $$\reg(R/(J+\wp^a)) \le \delta.$$
\end{lemma}

The following lemma has been proved in \cite{Th6}.

\begin{lemma} \label{l31} \cite[Lemma 3.3]{Th6} Let $X=\{P_1,\ldots , P_s\}$
be a set of distinct points
in $\mathbb P^n$, and $m_1,\ldots , m_s$ be positive integers. Put $I =
\wp^{m_1}_1\cap \cdots\cap \wp^{m_s}_s$. If $Y=\{P_{i_1}, \ldots, P_{i_r}\}$ is
a subset of $X$ and $J = \wp^{m_{i_1}}_{i_1}\cap \cdots\cap
\wp^{m_{i_r}}_{i_r}$, then
$$\reg(R/I) \ge \reg(R/J).$$

\end{lemma}

\section{Lower bound for the regularity index of fat points}

Let $X=\{P_1,\ldots , P_s\}$ be a set of distinct points in $\mathbb P^n$ and
$m_1,\ldots , m_s$  be positive integers. Let $n_1,\ldots , n_s$ be non-negative
integers with $(n_1, \ldots, n_s) \ne (0, \ldots, 0)$ and $m_i\ge n_i$ for $i=1,
\ldots, s$. Put $I = \wp^{m_1}_1\cap \cdots\cap \wp^{m_s}_s$, $N =
\wp^{n_1}_1\cap \cdots\cap \wp^{n_s}_s$  ($\wp_i^{n_i}=R$ if $n_i=0$). Then we
have $$e(R/I) \ge e(R/N) \text{ and } H_{R/I}(t) \ge H_{R/N}(t).$$ So, we can
not compare $\reg(R/I)$ with $\reg(R/N)$ by definition of the regularity index
of fat points. In Proposition \ref{p33}  we will prove that $\reg(R/I) \ge
\reg(R/N)$.

The first, we get the following result.

\begin{lemma} \label{l32} Let $X=\{P_1,\ldots , P_s\}$ be a set of distinct points
in $\mathbb P^n$ and $m_1,\ldots , m_s$, $n_1,\ldots , n_s$ be positive integers
with $m_i\ge n_i$ for $i=1, \ldots, s$. Put $I = \wp^{m_1}_1\cap \cdots\cap
\wp^{m_s}_s$ and $N = \wp^{n_1}_1\cap \cdots\cap \wp^{n_s}_s$, then
$$\reg(R/I) \ge \reg(R/N).$$
\end{lemma}

\begin{proof} In case $m_i= n_i$ for $i=1, \ldots, s$, we have the equality.
In case there exists $j$ such that $m_j > n_j$, we may assume that $m_s > n_s$.
Put $I_1=\wp^{m_1}_1\cap \cdots\cap \wp^{m_{s-1}}_{s-1} \cap
\wp^{m_s-1}_s$. We will prove $\reg(R/I)\ge \reg(R/I_1)$.

Put $J=\wp^{m_1}_1\cap \cdots\cap \wp^{m_{s-1}}_{s-1}$. By Lemma \ref{l22} we
have
\begin{align*} &\reg(R/(J+\wp^{m_s}_s))\le b \\
 &\Leftrightarrow X_0^{b-i}M \in J+\wp_s^{i+1} \text{ for every }
 M=X_1^{c_1}\cdots X_n^{c_n}, c_1+\cdots+c_n=i, i=1, \ldots, m_s-1\\
 &\Rightarrow X_0^{b-i}M \in J+\wp_s^{i+1} \text{ for every }
 M=X_1^{c_1}\cdots X_n^{c_n}, c_1+\cdots+c_n=i
 , i=1, \ldots, m_{s}-2\\
 &\Leftrightarrow \reg(R/(J+\wp^{m_s-1}_s))\le b.
\end{align*}
This implies $\reg(R/(J+\wp^{m_s}_s)) \ge \reg(R/(J+\wp^{m_s-1}_s))$. By Lemma
\ref{l21} we have \begin{align*} \reg(R/I)&=\max\left\{m_{s}-1, \reg(R/J),
\reg(R/(J+\wp^{m_s}_s)) \right\},\\
\reg(R/I_1)&=\max\left\{m_{s}-2, \reg(R/J), \reg(R/(J+\wp^{m_s-1}_s)) \right\}.
\end{align*}
Therefore, we get $$\reg(R/I)\ge \reg(R/I_1).$$ By induction on $m_s$ we get
$$\reg(R/I) \ge \reg(R/(\wp_1^{m_1} \cap \cdots \cap \wp_{s-1}^{m_{s-1}}\cap
\wp_s^{n_s})).$$ By induction on number of points we get
$$\reg(R/I) \ge \reg(R/N).$$

\end{proof}

From the above lemma and Lemma \ref{l31} we get the following proposition.

\begin{proposition}\label{p33} Let $X=\{P_1,\ldots , P_s\}$ be a set of distinct points
in $\mathbb P^n$ and $m_1,\ldots , m_s$  be positive integers. Let $n_1,\ldots ,
n_s$ be non-negative integers with $(n_1, \ldots, n_s) \ne (0, \ldots, 0)$ and
$m_i\ge n_i$ for $i=1, \ldots, s$. Put $I = \wp^{m_1}_1\cap \cdots\cap
\wp^{m_s}_s$, $N = \wp^{n_1}_1\cap \cdots\cap \wp^{n_s}_s$  ($\wp_i^{n_i}=R$ if
$n_i=0$). We have
$$\reg(R/I) \ge \reg(R/N).$$
\end{proposition}

\medskip

Now by using the above proposition we show a lower bound for the regularity
index of fat points. We recall that all rational normal curves in $\mathbb P^j$
are isomorphic under a linear change of coordinates. Hence, we may assume that
the points are on the rational normal curve in $\mathbb P^j$ with the parametric
equations
$$X_0=t^j, X_1=t^{j-1}u, \ldots, X_{j-1}=tu^{j-1}, X_j=u^j.$$

Let $Q_1,\ldots, Q_r$ be distinct points in $\mathbb P^n$. If there exist a
linear change of coordinates $\varphi$ of $\mathbb P^n$ such that the
coordinates of the points $\varphi(Q_1), \ldots, \varphi(Q_r)$ satisfy the
parametric equations
$$X_0=t^j, X_1=t^{j-1}u, \ldots, X_{j-1}=tu^{j-1}, X_j=u^j, X_{j+1}=0, \ldots,
X_n=0,$$ then we said that $Q_1, \ldots, Q_r$ are in Rnc-j in $\mathbb P^n$. So,
without loss of generality, we may say that if $Q_1, \ldots, Q_r$ are in Rnc-j
in $\mathbb P^n$, then their coordinates satisfy the above parametric equations.

The points lying on a line are Rnc-1. The points lying on a rational normal
curve in $\mathbb P^n$ are Rnc-n in $\mathbb P^n$.

\par\medskip

Let $Z=m_1P_1+\cdots +m_{s}P_{s}$ be a set of fat points in $\mathbb P^n$. Then
the set
$$\{P_{i_1}, \ldots, P_{i_q}\in \{P_1, \ldots, P_s\}| P_{i_1}, \ldots, P_{i_q}\text{ are in }
Rnc-1\}$$ is non-empty.

\medskip

If $\frac mj$ is a rational number, we denote by $[\frac mj]$ it's integer part.

The following theorem shows a lower bound for the regularity index of fat
points.

\begin{theorem}\label{tlb}  Let $Z=m_1P_1+\cdots +m_{s}P_{s}$ be a set of fat
points in $\mathbb P^n$. Then, $$\reg(Z) \ge \max \{ D_j |\  j=1,\ldots, n \},$$
where $$D_j = \max\left\{\left[\frac{\sum_{l=1}^q m_{i_l}+ j- 2}{j}\right] |\
P_{i_1}, \ldots , P_{i_q} \text{ are in Rnc-j}\right\}.$$
\end{theorem}

\begin{proof}
Suppose that points $P_{i_1}, \ldots, P_{i_q}$ of $\{P_1, \ldots, P_n\}$ are
Rnc-j in $\mathbb P^n$. We may assume that $m_{i_1}\ge \cdots \ge m_{i_q}$
(after relabeling the points, if necessary). Let $\wp_{i_1}, \ldots, \wp_{i_q}$
be the homogeneous prime ideals of $R$ corresponding to the points $P_{i_1},
\ldots, P_{i_q}$. Put
$$J=\wp_1^{m_{i_1}}\cap \cdots \cap \wp_{i_{q}}^{m_{i_{q}}}.$$
By Lemma \ref{l31} we have
$$\reg(Z)\ge \reg(R/J).$$ We will prove that
$$\reg(R/J)\ge \left[\frac{\sum_{l=1}^q m_{i_l}+ j-
2}{j}\right].$$

Since the points $P_{i_1}, \ldots, P_{i_q}$ are in Rnc-j in $\mathbb P^n$, we
may assume that their coordinates satisfying parametric equations:
$$X_0=t^j, X_1=t^{j-1}u, \ldots, X_{j-1}=tu^{j-1}, X_j=u^j, X_{j+1}=0, \ldots,
X_n=0$$ and the points $P_{i_q}=(1, 0, \ldots, 0)$.  Then $\wp_{i_q}=(X_1,
\ldots, X_n)$. Put
$$J_1=\wp_{i_1}^{m_{i_1}}\cap \cdots \cap \wp_{i_{q-1}}^{m_{i_{q-1}}}.$$
The first, we will prove that
$$\reg(R/(J_1+\wp_{i_q}^{m_{i_q}})) \ge \left[\frac{\sum_{l=1}^q m_{i_l}+ j-
2}{j}\right].$$ Put $T=\left[\frac{\sum_{l=1}^q m_{i_l}+ j- 2}{j}\right]$.
Consider the monomial $X_0^{T-m_{i_q}}X_1^{m_{i_q}-1}$. If
$$X_0^{T-m_{i_q}}X_1^{m_{i_q}-1} \in J_1 + \wp_{i_q}^{m_{i_q}},$$ then there
exists a form $h \in \wp_{i_q}^{m_{i_q}}$ of degree $T-1$ such that
$$X_0^{T-m_{i_q}}X_1^{m_{i_q}-1} + h \in J_1.$$ Since $X_0^{T-m_{i_q}}X_1^{m_{i_q}-1} \in
\wp_{i_q}^{m_{i_q}-1}$ and $h \in \wp_{i_q}^{m_{i_q}} \subset
\wp_{i_q}^{m_{i_q}-1}$, we have
$$X_0^{T-m_{i_q}}X_1^{m_{i_q}-1} + h \in J_1\cap \wp_{i_q}^{m_{i_q}-1}=\wp_{i_1}^{m_{i_1}}
\cap \cdots \cap \wp_{i_{q-1}}^{m_{i_{q-1}}} \cap \wp_{i_q}^{m_{i_q}-1}.$$ Since
$m_{i_1}+\cdots +m_{i_q} -1> j(T-1)$, by Bezout's theorem we have
$$X_0^{T-m_{i_q}}X_1^{m_{i_q}-1} + h$$ vanishing on the points $(1, \lambda,
\ldots, \lambda^j, 0, \ldots, 0)\in \mathbb P^n$, for every $\lambda$ in the
field $k$. This implies
$$\lambda^{m_{i_q}-1}+h(1, \lambda, \ldots, \lambda^j, 0, \ldots, 0)=0$$ for every $\lambda\in k$.
Since $h\in \wp_{i_q}^{m_{i_q}}=(X_1, \ldots, X_n)^{m_{i_q}}$, we have $h(1,
\lambda, \ldots, \lambda^j, 0, \ldots, 0)=0$ or  $h(1, \lambda, \ldots,
\lambda^j, 0, \ldots, 0)=\lambda^{m_{i_q}}g(\lambda)$, for some non-zero
polynomial $g\in k[x]$. Hence, $\lambda^{m_{i_q}-1}=0$ or
$\lambda^{m_{i_q}-1}+\lambda^{m_{i_q}}g(\lambda)=0$ for every $\lambda \in k$, a
contradiction. Thus, we get
$$X_0^{T-m_{i_q}}X_1^{m_{i_q}-1} \notin J_1 + \wp_{i_q}^{m_{i_q}}.$$ By
Lemma \ref{l22} we have $$\reg (R/(J_1+\wp_{i_q}^{m_{i_q}})) \ge T.$$ Next, by
Lemma \ref{l21} we get \begin{align*} \reg(R/J)&=\max\{m_{i_q}-1, \reg(R/J_1),
\reg(R/(J_1+\wp^{m_{i_q}}_{i_q}))\}\\ &\ge \reg (R/(J_1+\wp_{i_q}^{m_{i_q}}))
\ge T.\end{align*}

The proof of Theorem \ref{tlb} is now completed. \end{proof}

\section{Application of lower bound}

The first, by using the lower bound we can compute the regularity index of fat
points whose support on a line. This formula was showed by E.D. Davis and A.V.
Geramita \cite[Corollary 2.3]{DG} by using another method.

\begin{proposition}\label{p41} Let $Z=m_1P_1+\cdots +m_{s}P_{s}$ be a set of fat
points in $\mathbb P^n$. If $P_1, \ldots, P_s$ lie on a line, then $$\reg(Z)=
m_1+\cdots+m_s-1.$$
\end{proposition}

\begin{proof} We may assume that $m_1 \ge \cdots \ge m_s$.
Since $P_1, \ldots, P_s$ lie on a line, we have $D_1=m_1+\cdots+m_s-1$. Put
$I=\wp_1^{m_1} \cap\cdots\cap \wp_s^{m_s}$. Then $\reg(Z)=\reg(R/I)$. We will
prove that
$$\reg(R/I)=D_1.$$ By Theorem \ref{tlb} we have
$$\reg(R/I) \ge \max \{ D_j |\
j=1, \ldots, n \}.$$ So, it suffices to prove by induction on $s$ that
$$\reg(R/I)\le D_1.$$ Put $J=\wp_1^{m_1} \cap\cdots\cap \wp_{s-1}^{m_{s-1}}$, by
the inductive assumption, we get
$$\reg(R/J)\le m_1+\cdots+m_{s-1}-1 \le D_1.\quad\quad\quad\quad (1)$$

 Choose $P_s=(1, 0, \ldots, 0)$,
then $\wp_s=(X_1, \ldots, X_n)$.  For $j=1, \ldots, s-1$, since $P_1, \ldots,
P_s$ lie on a line, there exists hyperplane, say $H_j$, passing throught $P_j$
and avoiding $P_s$. This implies
$$H_1^{m_1} \cdots H_{s-1}^{m_{s-1}} \in J.$$
Therefore, for every monomial $M=X_1^{c_1}X_2^{c_2} \cdots X_n^{c_n}$ of degree
$i$, $i=0, \ldots, m_s-1$, we have $$H_1^{m_1} \cdots H_{s-1}^{m_{s-1}}M \in
J.$$ Since $D_1 \ge \max\{m_1+\cdots +m_{s-1} +i | i=1, \ldots, m_s-1\}$, by
Lemma \ref{l23} we get
$$\reg(R/(J+\wp_s^{m_s})) \le D_1.\quad\quad\quad\quad (2)$$
From (1), (2) and Lemma \ref{l21} we get
$$\reg(R/I) \le D_1.$$

\end{proof}

Next, we consider the fat points whose support lie on two separate lines.

\begin{proposition}\label{p42} Let $Z=m_1P_1+\cdots +m_{s}P_{s}$ be a set of fat
points in $\mathbb P^n$. Assume that there exist two lines $l_1$ and $l_2$, $l_1
\cap l_2=\emptyset$ such that $P_1, \ldots, P_s \in l_1 \cup l_2$. Then
$$\reg(Z)= D_1,$$
where $D_1=\max\left\{\sum_{l=1}^q m_{i_l}-1 |\ P_{i_1}, \ldots , P_{i_q} \text{
lie on a line}\right\}.$
\end{proposition}

\begin{proof}  Put
$I=\wp_1^{m_1} \cap\cdots\cap \wp_s^{m_s}$. Then $\reg(Z)=\reg(R/I)$. We will
prove that
$$\reg(R/I)=D_1.$$ By Theorem \ref{tlb} we have
$$\reg(R/I) \ge D_1.$$ So, it suffices to prove that
$$\reg(R/I)\le D_1.$$ We may assume that $P_1, \ldots,P_r \in l_1$ and $P_{r+1},
\ldots,P_s \in l_2$. Choose $P_s=(1, 0, \ldots, 0)$, then $\wp_s=(X_1, \ldots,
X_n)$. Put $J=\wp_1^{m_1} \cap\cdots\cap \wp_{s-1}^{m_{s-1}}$ and
$a=\max\{m_j|j=1, \ldots, r\}$.  We consider two following cases.

\smallskip

\noindent {\it Case 1:} $r=s-1$. Then $P_1, \ldots, P_{s-1} \in l_1$. Let $H$ be
the hyperplane containing $l_1$ and avoiding $P_s$. For every monomial $M$ of
degree $i$ in $X_1, \ldots, X_n$, $i=1, \ldots, m_s-1$, we have $$H^a M \in J.$$
Since $D_1 \ge a+ m_s-1 \ge \max\{a+i| i=1, \ldots, m_s-1\}$, by Lemma \ref{l23}
we get
$$\reg(R/(J+\wp_s^{m_s})) \le D_1.$$
Since $P_1, \ldots, P_{s-1}$ lie on a line, we get $$\reg(R/J)=m_1+\cdots
+m_{s-1}-1 \le D_1.$$ Therefore, by Lemma \ref{l21} $$\reg(R/I)=\max\{m_s-1,
\reg(R/J), \reg(R/(J+\wp_s^{m_s}))\} \le D_1.$$

\smallskip

\noindent {\it Case 2:} $r\le s-2$. We may argue by induction on $s$. Put
$b=\max\{a, m_{r+1} +\cdots +m_{s-1}\}$. Let $H_j$ be the hyperplane containing
$P_{r+j}$, $l_1$ and avoiding $P_s$, $j=1, \ldots, s-r-1$. For every monomial
$M$ of degree $i$ in $X_1, \ldots, X_n$, $i=1, \ldots, m_s-1$, we have
$$H_1^{b-(m_{r+1} +\cdots +m_{s-1})}H_1^{m_{r+1}} \cdots H_{s-r-1}^{m_{s-1}}M \in J,$$
($H_1^0 =R$). Since $D_1 \ge \max\{b+i | i=1, \ldots, m_s-1\}$, by Lemma
\ref{l23} we get
$$\reg(R/(J+\wp_s^{m_s})) \le D_1.$$
By the inductive assumption, we get
$$\reg(R/J) \le D_1.$$
Therefore, by Lemma \ref{l21} $$\reg(R/I)=\max\{m_s-1, \reg(R/J),
\reg(R/(J+\wp_s^{m_s}))\} \le D_1.$$

\end{proof}

Now we consider a set of fat points whose support is in Rnc-j.

\begin{lemma}\label{l42} Let $Z=m_1P_1+\cdots +m_sP_s$ be a set of fat
points in $\mathbb P^n$. Suppose that $j$ is the least integer such that $P_1,
\ldots, P_s$ are in Rnc-j. If $t$ is an integer such that
$$t\ge \max\left\{m_l, \left[\frac{m_1+\cdots+m_{s-1}+j- 1}{j}\right] |\
l=1, \ldots, s-1 \right\},$$ then we can find $t$ hyperplanes, say $H_1, \ldots,
H_t$, avoiding $P_s$ sucth that
$$H_1\cdots H_t \in \wp_1^{m_1}\cap \cdots \cap \wp_{s-1}^{m_{s-1}}.$$

\end{lemma}

\begin{proof} Since the points $P_1, \ldots, P_s$ are in Rnc-j in $\mathbb P^n$,
we may assume that their coordinators satisfying parametric equations:
$$X_0=v^j, X_1=v^{j-1}u, \ldots, X_{j-1}=vu^{j-1}, X_t=u^j, X_{j+1}=0, \ldots,
X_n=0.$$ If $j=1$, then $P_1, \ldots, P_s$ lie on a line. For $j=1, \ldots,
s-1$, there exists a hyperplane, say $H_j$, passing throught $P_j$ and avoiding
$P_s$ . Then we have $t=m_1+\cdots+m_{s-1}$ hyperplanes
$\underset{m_1}{\underbrace{H_1, \ldots, H_1}}, \ldots,
\underset{m_{s-1}}{\underbrace{H_{s-1}, \ldots, H_{s-1}}}$ avoiding $P_s$ such
that
$$H_1^{m_1} \cdots H_{s-1}^{m_{s-1}} \in \wp_1^{m_1}\cap \cdots \cap \wp_{s-1}^{m_{s-1}}.$$
If $t\ge 2$,  then no $l+2$ points of $\{P_1, \ldots, P_s\}$ are on a $l$-plane
for $l<j$. This implies that there does not exist any $(j-1)$-plane containing
$j+1$ points of $\{P_1, \ldots, P_s\}$. We will prove the lemma by induction on
$\sum_{i=1}^{s-1}m_i$.

We may assume that $m_1 \ge \cdots\ge m_{s-1}$. Since  $j$ is the least integer
such that $P_1, \ldots, P_s$ are in Rnc-j, we have $j \le s-1$. Let $\sigma_1$
be the $(j-1)$-plane passing throught $P_1, \ldots, P_{j}$. Then $\sigma_1$
avoids $P_s$. Therefore, there is a hyperplane, say $L_1$, containing $\sigma_1$
and avoiding $P_s$.

\noindent Case $s-1=j$: Then $$L_1^t \in  \wp_1^{m_1}\cap \cdots \cap
\wp_{s-1}^{m_{1}} \subset \wp_1^{m_1}\cap \cdots \cap \wp_{s-1}^{m_{s-1}}.$$

\noindent Case $s-1\ge j+1$: Since $t \ge \left[\frac{m_1+\cdots+m_{s-1}+j-
1}{j}\right]$ and $m_1\ge \cdots \ge m_{s-1}$, we have
\begin{align*} t-1 &\ge \left[\frac{m_1+\cdots+m_{s-1}+j- 1}{j}\right]-1
\ge \left[\frac{(j+1)m_{j+1}- 1}{j}\right]\\
&\ge m_{j+1}. \end{align*} On the other hand, since $t \ge
\left[\frac{m_1+\cdots+m_{s-1}+j- 1}{j}\right]$, we get $$t-1 \ge
\left[\frac{(m_1-1)+\cdots+(m_j-1)+m_{j+1}+\cdots+m_{s-1}+j- 1}{j}\right].$$
Consider
$$Z_1=(m_1-1)P_1+\cdots+(m_j-1)P_{j}+m_{j+1}P_{j+1}+\cdots+m_{s-1}P_{s-1}+m_sP_s.$$
By the inductive assumption we can find $(t-1)$ hyperplanes, say $L_2, \ldots,
L_t$, avoiding $P_s$ such that
$$L_2 \cdots L_t \in \wp_1^{m_1-1}\cap \cdots\cap \wp_{j}^{m_j-1}\cap
\wp_{j+1}^{m_{j+1}}\cap \cdots \cap \wp_{s-1}^{m_{s-1}}.$$ Moreover, since
$L_1\in \wp_1 \cap \cdots\cap \wp_j$, we get
$$L_1L_2\cdots L_t \in \wp_1^{m_1}\cap \cdots \cap \wp_{s-1}^{m_{s-1}}.$$

\end{proof}

We can compute the regularity index of fat points whose support is in Rnc-j.

\begin{proposition}\label{p43} Let $Z=m_1P_1+\cdots +m_{s}P_{s}$ be a set of fat
points in $\mathbb P^n$. If $P_1, \ldots, P_s$ are in Rnc-t, then
$$\reg(Z) =\max \{D_j| j=1, \ldots, t\},$$ where $$D_j
= \max\left\{\left[\frac{\sum_{l=1}^q m_{i_l}+ j- 2}{j}\right] |\ P_{i_1},
\ldots , P_{i_q} \text{ are in Rnc-j}\right\}.$$\end{proposition}

\begin{proof} We may assume that $m_1 \ge \cdots \ge m_s$.  We will argue by
induction on $s$. If $s=1$, then $\reg(Z)=m_1-1=D_1$. If $s\ge 2$, we consider
two following cases:

\noindent Case $t=1$: Then $P_1, \ldots, P_s$ lie on a line and
$D_1=m_1+\cdots+m_s-1=\max \{ D_j |\ j=1, \ldots, n \}$. By Proposition
\ref{p41} we have $\reg(Z)=D_1$.

\noindent Case $t\ge 2$: Since $P_1, \ldots, P_s$ are in Rnc-t, there is the
least integer $p \le t$ such that $P_1, \ldots, P_s$ are in Rnc-p. Then
\begin{align*} D_1&=m_1+m_2-1 \ge D_2 \ge \cdots \ge D_{p-1},\\
D_p&=\left[\frac{m_1+\cdots+m_{s}+ p- 2}{p}\right] \ge D_{p+1} \ge \cdots \ge
D_n.\end{align*} So, $\max \{ D_j |\ j=1, \ldots, n \}=\max \{ D_j |\ j=1,
\ldots, t \}=\max\{D_1, D_p\}$. Hence, by Theorem \ref{tlb} we get
$$\reg(Z) \ge \max\{D_1, D_p\}.$$
It sufficies to prove that $$\reg(Z) \le \max\{D_1, D_p\}.$$  Put
$Z_1=m_1P_1+\cdots+m_{s-1}P_{s-1}$ and $J=\wp_1^{m_1} \cap \cdots \cap
\wp_{s-1}^{m_{s-1}}$. We have $\reg(Z_1)=\reg(R/J)$. By inductive hypothesis we
have
$$\reg(Z_1)=\max \{D'_j| j=1, \ldots, t\},$$ where $$D'_j
= \max\left\{\left[\frac{\sum_{l=1}^q m_{i_l}+ j- 2}{j}\right] |\ P_{i_1},
\ldots , P_{i_q}\in \{P_1, \ldots, P_{s-1}\}  \text{ and are in
Rnc-j}\right\}.$$ Since $\{P_1, \ldots, P_{s-1}\} \subset \{P_1, \ldots,
P_{s-1}, P_s\}$, we have $D'_j \le D_j$ for $j=1, \ldots, t$. Therefore, we get
$$\reg(R/J) \le \max\{D_1, D_p\}. \quad\quad\quad\quad (3)$$
Consider $R/(J+\wp_s^{m_s})$. We may assume that $P_s=(1, 0, \ldots, 0)$,
$P_1=(0, \underset{2}{\underbrace{1}}, 0, \ldots, 0)$, $\ldots$, $P_p=(0,
\ldots, 0, \underset{p+1}{\underbrace{1}}, 0, \ldots, 0)$. For every monomial
$M=X_1^{c_1}\cdots X_n^{c_n}$, $c_1+\cdots+c_n=i$, $i = 0,\ldots , m_s-1$. Put
\[m'_l=\begin{cases}
m_l-i+c_l &\mbox{ for } l=1, \ldots, p,\\
m_l &\mbox{ for } l=p+1, \ldots, s-1.
  \end{cases}\]
Put $J'=\wp_1^{m'_1} \cap \cdots \cap \wp_{s-1}^{m'_{s-1}}$. By Proposition
\ref{l31} we can find
$$t= \max\left\{m'_l, \left[\frac{m'_1+\cdots+m'_{s-1}+p-
1}{p}\right] | l=1, \ldots, s-1\right\}$$ hyperplanes, say $H_1, \ldots, H_t$,
avoiding $P_s$ such that
$$H_1 \cdots H_t \in J'.$$ Since $M\in \wp_1^{i-c_1}
\cap \cdots \cap \wp_p^{i-c_p}$ and $J'=\wp_1^{m_1-i+c_1} \cap \cdots\cap
\wp_p^{m_p-i+c_p}\cap \wp_{p+1}^{m_{p+1}} \cdots \cap \wp_{s-1}^{m_{s-1}}$, we
get
$$H_1 \cdots H_t M \in J.$$
By Lemma \ref{l23} we get
$$\reg(R/(J+\wp_s^{m_s})) \le \max\{t+i|i=1, \ldots, m_s-1\} \le \max\{D_1,
D_p\}. \quad \quad \quad (4)$$ Put $I=J \cap \wp_s^{m_s}$. We have
$\reg(Z)=\reg(R/I)$. From (3), (4) and Lemma \ref{l21} we have $$\reg(Z) \le
\max\{D_1, D_p\}.$$

\end{proof}

In Proposition \ref{p43}, if $P_1, \ldots, P_s$ are in Rnc-n in $\mathbb P^n$,
then we get Proposition 7 in \cite{CTV}.

\par\noindent
Phan Van Thien, \\ Department of Mathematics, Hue Normal University,  \\ Vietnam
\\ E-mail: tphanvannl@yahoo.com

\end{document}